\documentclass[11pt,english]{amsart}
\usepackage[T1]{fontenc}
\usepackage[latin1]{inputenc}
\usepackage{amstext}
\usepackage{amsthm}
\usepackage{amssymb}

\makeatletter
\numberwithin{equation}{section}
\numberwithin{figure}{section}
\theoremstyle{plain}
\newtheorem{thm}{\protect\theoremname}
\theoremstyle{remark}
\newtheorem{rem}[thm]{\protect\remarkname}
\theoremstyle{plain}
\newtheorem{lem}[thm]{\protect\lemmaname}

\usepackage{tikz}
\usepackage{pgfplots}

\makeatother

\makeatother

\usepackage{babel}
\providecommand{\lemmaname}{Lemma}
\providecommand{\remarkname}{Remark}
\providecommand{\theoremname}{Theorem}

\begin{document}


\addtolength{\textwidth}{0mm}
\addtolength{\hoffset}{-0mm} 
\addtolength{\textheight}{0mm}
\addtolength{\voffset}{-0mm} 


\global\long\def\CC{\mathbb{C}}%
 
\global\long\def\BB{\mathbb{B}}%
 
\global\long\def\PP{\mathbb{P}}%
 
\global\long\def\QQ{\mathbb{Q}}%
 
\global\long\def\RR{\mathbb{R}}%
 
\global\long\def\FF{\mathbb{F}}%

\global\long\def\DD{\mathbb{D}}%
 
\global\long\def\NN{\mathbb{N}}%
\global\long\def\ZZ{\mathbb{Z}}%
 
\global\long\def\HH{\mathbb{H}}%
 
\global\long\def\Gal{{\rm Gal}}%
\global\long\def\OO{\mathcal{O}}%
\global\long\def\pP{\mathfrak{p}}%

\global\long\def\pPP{\mathfrak{P}}%
 
\global\long\def\qQ{\mathfrak{q}}%
 
\global\long\def\mm{\mathcal{M}}%
 
\global\long\def\aaa{\mathfrak{a}}%
 
\global\long\def\a{\alpha}%
 
\global\long\def\b{\beta}%
 
\global\long\def\d{\delta}%
 
\global\long\def\D{\Delta}%
 
\global\long\def\L{\Lambda}%
 
\global\long\def\g{\gamma}%

\global\long\def\G{\Gamma}%
 
\global\long\def\d{\delta}%
 
\global\long\def\D{\Delta}%
 
\global\long\def\e{\varepsilon}%
 
\global\long\def\k{\kappa}%
 
\global\long\def\l{\lambda}%
 
\global\long\def\m{\mu}%
 
\global\long\def\o{\omega}%
 
\global\long\def\p{\pi}%
 
\global\long\def\P{\Pi}%
 
\global\long\def\s{\sigma}%

\global\long\def\S{\Sigma}%
 
\global\long\def\t{\theta}%
 
\global\long\def\T{\Theta}%
 
\global\long\def\f{\varphi}%
 
\global\long\def\deg{{\rm deg}}%
 
\global\long\def\det{{\rm det}}%

\global\long\def\Dem{Proof: }%
 
\global\long\def\ker{{\rm Ker\,}}%
 
\global\long\def\im{{\rm Im\,}}%
 
\global\long\def\rk{{\rm rk\,}}%
 
\global\long\def\car{{\rm car}}%
\global\long\def\fix{{\rm Fix( }}%

\global\long\def\card{{\rm Card\  }}%
 
\global\long\def\codim{{\rm codim\,}}%
 
\global\long\def\coker{{\rm Coker\,}}%
 
\global\long\def\mod{{\rm mod }}%

\global\long\def\pgcd{{\rm pgcd}}%
 
\global\long\def\ppcm{{\rm ppcm}}%
 
\global\long\def\la{\langle}%
 
\global\long\def\ra{\rangle}%

\global\long\def\Alb{{\rm Alb(}}%
 
\global\long\def\Jac{{\rm Jac(}}%
 
\global\long\def\Disc{{\rm Disc(}}%
 
\global\long\def\Tr{{\rm Tr(}}%
 
\global\long\def\NS{{\rm NS(}}%
 
\global\long\def\Pic{{\rm Pic(}}%
 
\global\long\def\Pr{{\rm Pr}}%

\global\long\def\Km{{\rm Km(}}%
\global\long\def\rk{{\rm rk(}}%
\global\long\def\Hom{{\rm Hom(}}%
 
\global\long\def\End{{\rm End(}}%
 
\global\long\def\aut{{\rm Aut}}%
 
\global\long\def\SSm{{\rm S}}%
\global\long\def\psl{{\rm PSL}}%

\global\long\def\cu{{\rm (-2)}}%

\title{Conic configurations via dual of quartic curves}
\author{Xavier Roulleau}
\begin{abstract}
We construct special conic configurations from some point configurations
which are the singularities of the dual of a quartic curve. 
\end{abstract}

\maketitle

\section*{Introduction}

The aim of this paper is to report some experiments on the construction
of conic arrangements in the plane. Let us recall the very classical
(and elementary) fact that through $5$ points in general position
there is a unique irreducible conic and there is no conic through
$6$ points in general position. In this note, we study and construct
sets $\mathcal{C}$ of conics, each of which containing at least $6$
points of a special set $\mathcal{P}$ of points. The idea and basic
strategy for obtaining such set of points grew-up from the following
example: \\
Let $C=\{f=0\}\hookrightarrow\PP^{2}$ be a smooth cubic curve and
let $\check{C}\hookrightarrow\check{\PP^{2}}$ be the dual curve (the
image of $C$ by the gradient map $(\frac{df}{dx}:\frac{df}{dy}:\frac{df}{dz})$).
The singularities of the degree $6$ curve $\check{C}$ is a set $\mathcal{P}_{9}$
of $9$ cusps (corresponding to the $9$ flex lines of the cubic $C$).
It turns out that a set of irreducible conics $\mathcal{C}_{12}$
that contains at least $6$ points in $\mathcal{P}_{9}$ has cardinality
$12$, the sets $(\mathcal{P}_{9},\mathcal{C}_{12})$ of points and
conics form an interesting configuration of type 
\[
(9_{8},12_{6}),
\]
described in \cite{RS}. In that paper, that conic arrangement appears
in the context of generalized Kummer surfaces: the desingularization
of the double cover branched over $\check{C}$ is a K3 surface $X$
on which the pull-back of the conics are union of $(-2)$-curves which
have as well an interesting configuration. The freeness of the arrangement
of curves $\mathcal{C}_{12}$ is studied in \cite{PS}, where we learned
that this configuration has been also independently discovered in
\cite{DLPU}. 

That example of $12$ conics and $9$ points provides an interesting
example of points and conics in special position and motivated us
to find other conic configurations with such a property. The dual
curve $\check{C}$ of a smooth degree $>2$ curve $C$ is always a
curve with singularities, and from the above example one may expect
that the singular set is in special position with respect to the conics. 

In this paper, we study the conics that go through at least $6$ of
the singularities of the dual of the Fermat quartic, the Klein quartic
and a randomly chosen quartic curve which has $28$ rational bitangent
(corresponding to these $28$ bitangent, there are $28$ nodes on
the dual curve). In each cases, we find a high number of smooth conics
that contain at least $6$ points in the singularity set of the curve.
We compute the so-called Harbourne constant of the various conic arrangements
obtained. 

The literature on line arrangements is abundant, and their use occur
in many branches of mathematics. By example, special line arrangements
are the main ingredient for the construction of ball quotient surfaces
by Hirzebruch (see \cite{Hirzebruch}). One of the most fascinating
features of some line configurations and their set of singular points
is their high degrees of symmetries. In this note, we obtain as well
some very symmetric configurations of conics. 

Let us describe a bit more the results obtained on a example. Let
$\mathcal{P}_{28}$ be the set of $28$ rational points corresponding
dually to $28$ rational bitangent of a particular smooth quartic
curve $Q_{4/\QQ}$. The set $\mathcal{C}_{1008}$ of conics in special
position with respect to $\mathcal{P}_{28}$ has cardinality $1008$.
Each conic contains $6$ points in $\mathcal{P}_{28}$, each point
lies on $216$ conics, that gives a $(28_{216},1008_{6})$-configuration.
This configuration is very symmetric, by example it contains $63$
$(16_{6})$ and $63$ $(12_{16},32_{6})$ points and conics sub-configurations.
We found all these configurations for various examples of quartics
$Q_{4}$ ; we may ask if that is always the case for generic quartics
and what could be the meaning of that phenomenon relative to these
quartics.

One also may search for special configurations of higher degree curves,
the next example being irreducible smooth cubic curves containing
at least $10$ points of a given set of points in $\PP^{2}$. Searching
such configuration is much more time consuming for the computer, but
we give one example of a configuration of $63$ smooth cubic curves
defined over $\QQ$ and $28$ rational points in the plane (which
is in fact the above set $\mathcal{P}_{28}$) forming a $(28_{27},63_{12})$-configuration.
Apart from the $28$ points of multiplicity $27$, that arrangement
of $63$ cubics has $1008$ triple points and $4725$ double points. 

This note is organized as follows: in the first section we recall
some preliminary materials, in the second, we describe special conic
configurations obtained as the set of conics containing $6$ or more
points in the singular set of the dual of a quartic curve. In the
third section, we study configurations of conics that contain $6$
or more points in the singular set of some line arrangements, namely
the mirror of complex reflection groups. The last section is an appendix,
where we give the algorithms (for magma software V2.24) used to compute
the equations of the sets of conics, and the singularities of their
union. 

\textbf{Acknowledgements} The author is grateful to Piotr Pokora for
its comments on a first version of this paper, to the Referee for
a careful reading and suggestions, and to the Max Planck Institute
for Mathematics in Bonn for its hospitality and financial support.

\section{Preliminaries}

\subsection{Notations and conventions}

The index $n\in\NN$ of a set $\mathcal{S}_{n}$ indicates the cardinality
of the set. 

If $\mathcal{P}_{n}$ is a set of points in $\PP^{2}$ and $\mathcal{\ensuremath{C}}_{m}$
an arrangement of curves (by which we mean a set of plane curves)
such that each point $p\in\mathcal{P}$ is on $r$ curves in $\mathcal{C}_{m}$
and each curve $C\in\mathcal{C}_{m}$ contains $s$ points in $\mathcal{P}_{n}$,
we say that the sets $(\mathcal{P}_{n},\mathcal{C}_{m})$ form a $(n_{r},m_{s})$-configuration.
When $n=m$ (and thus $r=s$), one speaks of $n_{r}$-configurations.
For more informations and properties on abstract configurations in
algebraic geometry, we refer to \cite{Dolgachev}.

We remark that there is a natural identification between an arrangement
and the associated reduced divisor which is the union of the curves
in the arrangement; we will make no distinction between these two
notions.

For an arrangement $\mathcal{C}$ of smooth irreducible curves and
for $r>1$ an integer, $t_{r}=t_{r}(\mathcal{C})$ denotes the number
of points of multiplicity $r$ of the curve $\mathcal{C}$. That differs
a bit from Hirzebruch's notations \cite{Hirzebruch}, where $t_{r}$
stands for the number of ordinary singularities of multiplicity $r$.

Most of the time, we tried to use the following policy: when we study
all the singularities of $\mathcal{C}$, a union of curves, we use
the words ``curve arrangement'', when we study a pair $(\mathcal{P},\mathcal{C})$
where $\mathcal{C}$ is a union of curves, and $\mathcal{P}$ is a
fixed sub-set of the singular set, we use the words ``curve configuration''. 

By 'conic', we mean a degree $2$ smooth plane curve. 

\subsection{Harbourne constant of a curve}

The $H$-constants have been introduced in \cite{BRHH}. Let $C$
be a curve on a smooth surface $X$, let $f:\tilde{X}\to X$ be the
blow-up of $X$ at the singular points of $C$ (supposed non-empty,
of cardinality $s$) and let $\tilde{C}\hookrightarrow\tilde{X}$
be the strict transform of the curve $C$ by $f$. The Harbourne constant
of $C$ at its singular set is
\[
H(C)=\frac{(\tilde{C})^{2}}{s}.
\]
By example, if $\mathcal{C}=C_{1}+\dots+C_{k}$ is a degree $d$ plane
curve which is the union of smooth curves $C_{j}$, having for $r>1$,
$t_{r}$ points of multiplicity $r$, one has 
\[
H(C)=\frac{d^{2}-\sum t_{r}r^{2}}{\sum t_{r}}.
\]
If $\mathcal{C}_{k}$ is the union of $k$ degree $d$ smooth plane
curves in general position then $H(\mathcal{C}_{k})>-2$, and $\lim_{k}H(\mathcal{C}_{k})=-2$.
Obtaining arrangements of plane curves with $H$-constant lower than
$2$ is more difficult and it seems that no general construction is
known (see \cite{BRHH} for results on line arrangements). 

\section{Dual of quartic curves}

\subsection{The Fermat quartic}

Consider the curve $C=\{x^{4}+y^{4}+z^{4}=0\}$, its dual is the
curve

\[
\check{C}=\{x^{12}+3x^{8}y^{4}+3x^{4}y^{8}+y^{12}+3x^{8}z^{4}-21x^{4}y^{4}z^{4}+3y^{8}z^{4}+3x^{4}z^{8}+3y^{4}z^{8}+z^{12}=0\}.
\]
The set $\mathcal{P}_{28}$ of singular points of $\check{C}$ is
the unions of the set $\mathcal{P}_{16}$ of nodes and the set $\mathcal{P}_{12}$
of $E_{6}$-singularities, where $\mathcal{P}_{16}$ is:
\[
\begin{array}{c}
(-1:-1:1),(-1:1:1),(-1:-i:1),(-1:i:1),(1:-1:1),\\
(1:1:1),(1:-i:1),(1:i:1),(-i:-1:1),(-i:1:1),(-i:-i:1),\\
(-i:i:1),(i:-1:1),(i:1:1),(i:-i:1),(i:i:1).
\end{array}
\]
and $\mathcal{P}_{12}$ is:
\[
\begin{array}{c}
(0:-w:1),(0:w:1),(0:-w^{3}:1),(0:w^{3}:1),(-w:0:1),(w:0:1),\\
(-w^{3}:0:1),(w^{3}:0:1),(-w:1:0),(w:1:0),(-w^{3}:1:0),(w^{3}:1:0),
\end{array}
\]
here $i^{2}=-1$, $w=\frac{\sqrt{2}}{2}(1+i)$ (thus $w^{2}=i$). 

We compute that the set $\mathcal{C}$ of irreducible conics that
contains at least $6$ points in $\mathcal{P}_{28}$ has cardinality
$2736$, these conics are as follows:
\begin{itemize}
\item $2616$ conics contain exactly $6$ points in $\mathcal{P}_{28}$,
\item $96$ conics contain exactly $7$ points in $\mathcal{P}_{28}$,
\item $24$ conics contain exactly $8$ points in $\mathcal{P}_{28}$.
\end{itemize}
Through a point in $\mathcal{P}_{16}$ there are $546$ conics in
$\mathcal{C}$ and through a point in $\mathcal{P}_{12}$ there are
$652$ conics in $\mathcal{C}$.

\subsubsection{A $(16_{24},64_{6})$ point-conic configuration}

Since the behavior of the number of conics through nodes and $E_{6}$-singularities
is different, let us study the set of irreducible conics that contain
at least $6$ of the $16$ nodes of $\check{C}$. We obtain that there
is a set $\mathcal{C}_{64}$ of $64$ conics that contain at least
$6$ points in $\mathcal{P}_{16}$. Through each of the $16$ points
there are $24$ conics in $\mathcal{C}_{64}$, and each conic contains
exactly $6$ points among the $16$. The sets $\mathcal{P}_{16}$
and $\mathcal{C}_{64}$ of points and conics form a 
\[
(16_{24},64_{6})
\]
configuration. The configuration has singularities 
\[
t_{2}=2832,\,\,t_{3}=96,\,\,t_{4}=72,\,\,t_{24}=16.
\]
The $H$-constant is 
\[
H=-\frac{772}{377}\simeq-2.04.
\]

\subsubsection{A $(12_{28},56_{6})$ point-conic configuration}

If instead we consider the set $\mathcal{P}_{12}$ of $12$ points
which are $E_{6}$ singularities of $\check{C}$, we obtain that there
is a set $\mathcal{C}_{56}$ of $56$ conics such that each conic
contains at least $6$ points in $\mathcal{P}_{12}$. In fact, each
conic contains exactly $6$ points and through each points there are
$28$ conics, thus, we obtain a 
\[
(12_{28},56_{6})
\]
configuration of points and conics. We can check that when $6$ points
in $\mathcal{P}_{12}$ are on a conic in $\mathcal{C}_{56}$, then
there is a unique conic in $\mathcal{C}_{56}$ that contain the $6$
complementary points. On the union of the conics there are $44$ non-ordinary
singularities, which are tacnodes (i.e. the local intersection multiplicity
is $2$). One has
\[
t_{2}=1520,\,t_{28}=12,
\]
therefore the $H$-constant is:
\[
H=\frac{112^{2}-1520\cdot2^{2}-12\cdot28^{2}}{1532}=-\frac{736}{383}\simeq-1.92.
\]

\begin{rem}
For an arrangement of $k$ conics with ordinary singularities one
has
\[
4\left(\begin{array}{c}
k\\
2
\end{array}\right)=\sum_{r\geq2}\left(\begin{array}{c}
r\\
2
\end{array}\right)t_{r},
\]
but in that example, the left hand side equals $6160$, whereas the
right hand side equals to $6056$. This is due to the fact that some
singularities are not ordinary. 
\end{rem}

\subsubsection{A $(32_{3},24_{4})$ point-line configuration}

In the above example there are $44$ non-ordinary singularities;
$12$ of them are in the set $\mathcal{P}_{12}$. Let us denote by
$\mathcal{P}_{32}$ the complementary set. The set of lines $\mathcal{L}_{24}$
lines that contain strictly more that $2$ points in $\mathcal{P}_{32}$
has cardinality $24$, the union of these $24$ lines has equation
\[
\begin{array}{c}
x^{16}y^{8}+x^{12}y^{12}+x^{8}y^{16}+x^{16}y^{4}z^{4}+2x^{12}y^{8}z^{4}+\\
2x^{8}y^{12}z^{4}+x^{4}y^{16}z^{4}+x^{16}z^{8}+2x^{12}y^{4}z^{8}+3x^{8}y^{8}z^{8}\\
+2x^{4}y^{12}z^{8}+y^{16}z^{8}+x^{12}z^{12}+2x^{8}y^{4}z^{12}\\
+2x^{4}y^{8}z^{12}+y^{12}z^{12}+x^{8}z^{16}+x^{4}y^{4}z^{16}+y^{8}z^{16}=0.
\end{array}
\]
The multiplicities of the singularities of the union are 
\[
t_{2}=96,\,\,t_{3}=32,\,\,t_{8}=3.
\]
Each of the $24$ lines contains $13$ singularities. The sets $\mathcal{P}_{32}$
and $\mathcal{L}_{24}$ form a 
\[
(32_{3},24_{4})
\]
configuration of lines and points. Each line contains $8$ double
points, $4$ triple points and one octuple point. The $H$-constant
is 
\[
H=-\frac{288}{131}\simeq-2.19.
\]

\subsubsection{$(12_{8},24_{4})$ and $(16_{6},24_{4})$ point-conic configurations}

The set $\mathcal{C}_{24}$ of $24$ conics that contain $8$ points
in $\mathcal{P}_{28}$ and the set $\mathcal{P}_{12}\subset\mathcal{P}_{28}$
form a 
\[
(12_{8},24_{4})
\]
configuration. The sets $\mathcal{C}_{24}$ and $\mathcal{P}_{16}\subset\mathcal{P}_{28}$
form a 
\[
(16_{6},24_{4})
\]
configuration. The union of the $24$ conics has $292$ singular points,
one has:
\[
t_{2}=216,\,t_{4}=48,\,t_{6}=16,\,t_{8}=12.
\]
There are $24$ points where the intersection is not transverse, at
these points only two conics go through (thus $t_{2}=196+24$). The
$H$-constant is 
\[
H=\frac{48^{2}-216\cdot2^{2}-48\cdot4^{2}-16\cdot6^{2}-12\cdot8^{2}}{292}=-\frac{168}{73}\simeq-2.30.
\]

\subsection{The Klein quartic}

Let $C=\{x^{3}y+y^{3}z+z^{3}x=0\}$, its dual $\check{C}$ is a
degree $12$ curve with a set of $52$ singularities, which is the
union of \\
- a set $\mathcal{P}_{28}$ of $28$ nodes, corresponding to the $28$
bitangent, \\
- a set $\mathcal{P}_{24}$ of $24$ cusps, which correspond to flex
points on the quartic.\\
These $52$ points are defined over $\QQ(\zeta_{7})$. 

\subsubsection{A $(28_{6},21_{8})$ point-conic configuration}

 There are $8$ lines that contain at least $3$ points of $\mathcal{P}_{28}$.
There are $3129$ irreducible conics that contain at least $6$ points.
These conics contain $6,7$ or $8$ points. The set $\mathcal{C}_{21}$
of conics that contain $8$ points in $\mathcal{P}_{28}$ has cardinality
$21$. Through each points in $\mathcal{P}_{28}$ there are $6$ conics,
thus the sets $\mathcal{P}_{28},\mathcal{C}_{21}$ form a 
\[
(28_{6},21_{8})
\]
configuration. The singularities are 
\[
t_{2}=420,\,t_{6}=28,
\]
which gives 
\[
H=-\frac{33}{16}\simeq-2.06.
\]
The intersections are transverse. The log-Chern ratio of the arrangement
is 
\[
\frac{1161}{543}\simeq2.13.
\]

\subsubsection{A $(24_{49},147_{8})$ point-conic configuration}

The set $\mathcal{C}_{147}$ of conics that contain $8$ points in
$\mathcal{P}_{24}$ has cardinality $147$; the sets $\mathcal{P}_{24}$
and $\mathcal{C}_{147}$ form a 
\[
(24_{49},147_{8})
\]
configuration. We have 
\[
t_{2}=14280,\,\,t_{3}=42,\,\,t_{49}=24,
\]
thus the $H$-constant is
\[
H=\frac{(147\cdot2)^{2}-14280\cdot2^{2}-42\cdot3^{2}-24\cdot49^{2}}{14346}=-\frac{4781}{2391}\simeq-1.99.
\]
There are some non-ordinary singularities. 

The subset $\mathcal{C}_{49}$ of conics in $\mathcal{C}_{147}$ that
go through the same point in $\mathcal{P}_{49}$ has cardinality $49$
and the singularities of their union are as follows:
\[
t_{2}=1127,\,\,t_{14}=2,\,\,t_{15}=21,\,\,t_{49}=1,
\]
thus the $H$-constant is
\[
H=\frac{98^{2}-1127\cdot2^{2}-2\cdot14^{2}-21\cdot15^{2}-49^{2}}{1151}=-\frac{2422}{1151}\simeq-2.10.
\]
There are $2$ non-ordinary singularities (through which pass $14$
conics). 

\subsection{A quartic with $28$ rational bitangents}

\subsubsection{Shioda's construction of quartics with $28$ rational bitangents}

To a set $(u_{1},\dots,u_{6})$ of rational numbers in a specified
open set of $\mathbb{A}^{6}$, Shioda (see \cite{Shioda}) associates
a quartic with $28$ rational bitangents and gives the coefficients
of the equations of the bitangents as degree $4$ polynomial functions
of the $u_{i}$'s. 

The affirmations of the following section have been obtained for the
quartic $Q_{4}$ associated to the randomly chosen parameters:
\[
u_{1}=-\frac{10}{9},\,\,u_{2}=\frac{1}{9},\,\,u_{3}=-\frac{6}{23},\,\,u_{4}=-\frac{7}{9},\,\,u_{5}=\frac{10}{9},\,\,u_{6}=\frac{9}{19}.
\]

We tested several other examples of parameters $u_{i}$'s and obtained
the same abstract point-conic configurations $(28_{216},1008_{6})$
and point-cubic configurations $(28_{27},63_{12})$ we describe below.
It would be interesting to understand if that is always the case for
the generic parameters $u_{i}$, and why that phenomenon happens.

\subsubsection{A $(28_{216},1008_{6})$ point-conic configuration}

Let $\mathcal{P}_{28}$ be the set of the $28$ rational nodes of
the dual curve of $Q_{4}$, nodes which corresponds to the $28$ rational
bitangents. Let $\mathcal{C}_{1008}$ be the set of conics that contain
at least $6$ points in $\mathcal{P}_{28}$. Each conic in $\mathcal{C}_{1008}$
contains $6$ points in $\mathcal{P}_{28}$ and through each point
there are $216$ conics, thus, we obtain a 
\[
(28_{216},1008_{6})
\]
configuration of points and conics. 
\begin{rem}
The $28$ points in $\mathcal{P}_{28}$ are in general position with
respect to the lines: no line contains $3$ points from $\mathcal{P}_{28}$. 
\end{rem}

\subsubsection{Sixty three $(12_{16},32_{6})$ point-conic sub-configurations}

Let us fix one of the conic $C_{o}\in\mathcal{C}_{1008}$ (we checked
that the following is true for any such conic). Consider the $22$
points that are the complementary set of points to the $6$ points
in $\mathcal{P}_{28}$ on $C_{o}$. Then there are $152$ conics that
contain at least $6$ points in that set of $22$ points, moreover
there are $12$ points through which there are $51$ conics and $10$
points through which there are $30$ conics. Consider that set $\mathcal{P}_{12}=\mathcal{P}_{12}(C_{o})$
of $12$ points. Then there is a set of $32$ conics in $\mathcal{C}_{1008}$
such that each conic contains $6$ points in $\mathcal{P}_{12}$,
and through each point in $\mathcal{P}_{12}$ there are $16$ conics,
thus, we get a configuration 
\[
(12_{16},32_{6})
\]
of points and conics. 

We checked that on the union of the $32$ conics, the $12$ points
have multiplicity $16$, there are $544$ double points and no other
intersection points. The singularities are ordinary and the $H$-constant
is 
\[
H=\frac{64^{2}-544\cdot2^{2}-12\cdot16^{2}}{544+12}=-\frac{288}{139}\simeq-2.07.
\]
The configuration $(12_{16},32_{6})$ has the following symmetries:
the set of $16$ conics that go through one point $p_{0}$ contains
only $11$ points among the $12$. Let $q_{0}$ be the the complementary
point; then the set of $16$ conics that go through $q_{0}$ is complementary
to the set of $16$ conics through $p_{0}$. In that way, we get two
configurations of type $(10_{8},16_{5})$ for the points different
than $p_{0},q_{0}$. That new configuration has 
\[
t_{2}=80,\,t_{8}=10,\,t_{16}=1,
\]
thus 
\[
H=\frac{32^{2}-80\cdot2^{2}-10\cdot8^{2}-16^{2}}{80+10+1}=-\frac{192}{91}\simeq-2.109.
\]

\subsubsection{A $(28_{27},63_{12})$ point-cubic configuration}

It turns out that there are exactly $63$ sets of $12$ points of
the form $\mathcal{P}_{12}=\mathcal{P}_{12}(C_{o})$ (when $C_{o}$
varies in $\mathcal{C}_{1008}$). Moreover, we obtain that through
each set $\mathcal{P}_{12}$ of $12$ points there is a unique cubic
curve, which is smooth. Let us denote by $\mathcal{E}_{63}$ the set
of such cubic curves. The sets $(\mathcal{P}_{28},\mathcal{E}_{63})$
form a 
\[
(28_{27},63_{12})
\]
configuration. The $63$ elliptic curves in $\mathcal{E}_{63}$ are
not isomorphic, and by taking their reduction modulo a prime, one
can even show that they are not isogeneous. The singularities of the
cubic curve configuration $\mathcal{E}_{63}$ are ordinary with $t_{2}=4725,\,\,t_{3}=1008,\,\,t_{27}=28$
and $H$-constant: 
\[
H=-\frac{1809}{823}\simeq-2.19.
\]

\begin{rem}
The elliptic curves in $\mathcal{E}_{63}$ contain naturally $12$
rational points. These $12$ points are the origin, a $2$-torsion
point $t_{2}$, $5$ points $p_{1},\dots,p_{5}$ and the points $t_{2}+p_{j}$,
$1\leq j\leq5$. One may ask how large is its Mordell-Weil group;
computations with a computer never finished, but we made some experiments
which seems to show that the Mordell-Weil group has rank at least
$5$. 
\end{rem}

\subsubsection{Sixty three $(16_{6})$ point-conic configurations}

We defined above $63$ sets $\mathcal{P}_{12}$ of cardinality $12$.
For such fixed set $\mathcal{P}_{12}$, let us consider its complementary
$\mathcal{P}_{16}\subset\mathcal{P}_{28}$. We obtain that the set
$\mathcal{C}_{16}$ of conics in $\PP^{2}$ that contain (at least)
$6$ points in $\mathcal{P}_{16}$ has cardinality $16$ and the sets
$\mathcal{P}_{16},\mathcal{C}_{16}$ form a $16_{6}$ configuration:
each point is contained in $6$ conics, each conic contains $6$ points.
The singularities are ordinary, with
\[
t_{2}=240,\,t_{6}=16,
\]
its $H$-constant is $H=-2$. 

\section{Complex Reflection groups}

In the above conic configurations linked to the dual of a quartic,
the conics contain at most $8$ of the points in the special set $\mathcal{P}$.
This is why we also searched for conic configurations related to complex
reflection group: below we find a set of conics, each one containing
$10$ points of a set $\mathcal{P}$. We study exceptional three dimensional
irreducible complex reflection groups; there are $5$ such groups,
which are classically denoted by $G_{23},...,G_{27}$. 

\subsection{The Group $G_{23}$}

The complex reflection group $G_{23}$, of order $120$, has $15$
mirrors (of order $2$). The union of the $15$ mirrors has a set
$\mathcal{P}_{15}$ of $15$ double points, a set $\mathcal{P}_{10}$
of $10$ triple points and a set $\mathcal{P}_{6}$ of $6$ quintuple
points. There are $2345$ conics containing at least $6$ points in
the union $\mathcal{P}_{15}\cup\mathcal{P}_{10}\cup\mathcal{P}_{6}$. 

There are $25$ conics that contain at least $6$ points in $\mathcal{P}_{10}$.
That set of conics $\mathcal{C}_{25}$ and the set $\mathcal{\ensuremath{P}}_{10}$
form a $(10_{15},25_{6})$-configuration. The number and multiplicities
of the singularities are
\[
t_{2}=150,\,\,t_{15}=10
\]
and the $H$-constant is $-\frac{35}{16}\simeq-2.18$. 

\subsection{The reflection group $G_{24}$ of cardinality $168$, the automorphism
group of the Klein quartic}

The complex reflection group $G_{24}$ is the order $168$ automorphism
group of the Klein quartic curve. The union of the $21$ mirrors (of
order $2$) has a set $\mathcal{P}_{21}$ of $21$ quadruple points
and a set $\mathcal{P}_{28}$ of $28$ triple points. 

There are $133$ irreducible conics containing at least $6$ points
in $\mathcal{P}_{21}$. That set of conics is the union of the set
$\mathcal{C}_{21}$ of the $21$ conics that contain $8$ points in
$\mathcal{P}_{21}$ and the set $\mathcal{C}_{112}$ of conics that
contain $6$ points in $\mathcal{P}_{21}.$ The sets $\mathcal{P}_{21}$
and $\mathcal{C}_{21}$ form a
\[
21_{8}
\]
configuration, the number and multiplicities of the singularities
are
\[
t_{2}=168,\,\,t_{8}=21
\]
and the $H$-constant is $-\frac{4}{3}$. Some singularities are non-ordinary. 
\begin{rem}
The Klein configuration of $21$ lines and $21$ points is a $21_{4}$-configuration.
\end{rem}

The union of the $21$ lines and $21$ conics in $\mathcal{C}_{21}$
is a curve with 
\[
t_{2}=378,\,\,t_{3}=28,\,\,t_{12}=21
\]
and $H$-constant equal to $-117/61\simeq-1.91$. 

In \cite[Section 3]{PR} is given another arrangement of $21$ conics
related to the group $G_{24}$, however these $21$ conics are different
from ours; their union have singularities: $t_{2}=168$, $t_{3}=224$. 

\subsection{The reflection groups $G_{25}$ and $G_{26}$, the automorphism group
of the Hesse pencil}

The configuration of $21$ lines obtained as the mirrors (of order
$2$ and $3$) of the complex reflection group $G_{26}$ have singularity
set with 
\[
t_{2}=36,\,\,t_{4}=9,\,\,t_{5}=12.
\]
There are no irreducible conics containing at least $6$ points in
the set $\mathcal{P}_{9}$ of points with multiplicity $4$, neither
there are for the set $\mathcal{P}_{12}$ of $12$ points with multiplicity
$5$. However, through the set $\mathcal{P}_{21}=\mathcal{P}_{9}\cup\mathcal{P}_{12}$
there pass $108$ irreducible conics that contain $6$ points. We
denote by $\mathcal{C}_{108}$ that set of conics. Through each point
in $\mathcal{P}_{9}$ there are $24$ conics, and through each points
in $\mathcal{P}_{12}$ there are $36$ conics. Each of the conics
contains $2$ points in $\mathcal{P}_{9}$ and $4$ points in $\mathcal{P}_{12}$,
thus, we get 
\[
(9_{24},108_{2})\text{ and }(12_{36},108_{4})
\]
configurations. 

The group $G_{25}$ is implicitly studied since its projectivisation
is the same as $G_{26}$, and the $12$ mirrors of $G_{25}$ are the
$12$ order $3$ mirrors of $G_{26}$: the $21$ singularities of
the line arrangement of $G_{25}$ form the set $\mathcal{P}_{21}$.

\subsection{The Group $G_{27}$, the Valentiner-Wiman group }

The group $G_{27}$ has $45$ mirrors, the singularities of the
union of these $45$ lines are 
\[
t_{3}=120,\,\,t_{4}=45,\,\,t_{5}=36.
\]
Let $\mathcal{P}_{36}$ be the set of multiplicity $5$ singularities.
There are $13062$ irreducible conics that contain at least $6$ points
in $\mathcal{P}_{36}$. Through each point there are $2200$ such
conics. The conics contain $6,8$ or $10$ points in $\mathcal{P}_{36}$.
The set $\mathcal{C}_{72}$ of conics that contain $10$ points has
cardinality $72$. Through each point there are $20$ conics, thus
the sets $\mathcal{P}_{36}$, $\mathcal{C}_{72}$ form a 
\[
(36_{20},72_{10})
\]
configuration of points and conics. The singularities of that configuration
are 
\[
t_{2}=3312,\,\,t_{20}=36,
\]
there are $72$ tacnodes and the $H$-constant equals to $-\frac{64}{31}\simeq-2.06$.

The set of conics containing $8$ points has cardinality $270$, through
each points in $\mathcal{P}_{36}$ there are $60$ conics, thus, we
get a $(36_{60},270_{8})$ configuration. The remaining conics and
the points in $\mathcal{P}_{36}$ form a $(36_{2120},12720_{6})$-configuration.

\section{Appendix: algorithm used to compute the conics and their intersection
points}

Let us recall the following elementary facts:
\begin{lem}
Let us consider $5$ points in the plane. Assume that no three of
these points are on a line. Then there is a unique conic that contains
these points and that conic is irreducible.\\
If three of these points lie on a line, then the conic is reducible,
and may or may not be unique. If no four points are collinear, then
the five points define a unique conic (degenerate if three points
are collinear, but the other two points determine the unique other
line). If four points are collinear, then there is not a unique conic
passing through them -- one line passing through the four points,
and the remaining line passes through the other point, leaving $1$
parameter free. If all five points are collinear, then the remaining
line is free, which leaves $2$ parameters free.
\end{lem}

In the following, we describe two algorithms used to compute the set
of conics containing at least $6$ points of a set of points in $\PP^{2}$
and the singularity set of their union (with the multiplicities of
these points); we use Magma software V. 2.24. 

The entry of the first algorithm is a set $P_{0}$ of points in the
plane, defined over a field $K$. The output is the set of irreducible
reduced conics that contain at least $6$ points in $P_{0}$. This
algorithm is as follows:
\begin{verbatim}
function ConicsThrough6Points(Po)
  P2<x,y,z>:=Scheme(Po[1]);
  L2:=LinearSystem(P2,2);
\end{verbatim}
Here $L_{2}$ is the linear system of conics in the plane. Rather
than computing the set of conics containing $6$ points ($\#P_{0}$
choose $6$ possibilities), we start by studying the set of conics
that contains at least $5$ points in $P_{0}$ ($\#P_{0}$ choose
$5$ possibilities) and then we sort the conics that contain $6$
points. Doing so we gain some speed:
\begin{verbatim}
  E:=SetToSequence(Subsets({1..#Po}, 5));
  Conics5:=[LinearSystem(L2,[P2!Po[k] : k in q]) : q in E];
\end{verbatim}
The set Conics5 contains the linear systems of conics that goes through
the points $P_{0}[k_{1}],\dots,P_{0}[k_{5}]$, where $\{k_{1},\dots,k_{5}\}$
varies into the set $E$ of $5$-tuples of elements of $\{1,\dots,\#P_{o}\}$.
We then select among them the $1$-dimensional linear systems (a necessarily
condition in order that the conic is irreducible) and define the list
$Conics5P$ of conics that contains at least $5$ points:
\begin{verbatim}
  Conics5P:=[];
  for ru in [1..#E] do 
    q:=Sections(Conics5[ru]); 
    if #q eq 1 then 
      Append(~Conics5P,Conics5[ru]);  
    end if;
  end for;
  Conics5P:=[Scheme(P2,Sections(q)[1]): q in Conics5P}];
\end{verbatim}
Here we construct the set of irreducible reduced conics that contains
at least $6$ points in $P_{0}$:
\begin{verbatim}
  Conics6:={q : q in Conics5P | #{p : p in Po | p in q} gt 5}; 
  Conics6:=SetToSequence({q : q in Conics6 | IsIrreducible(q) 
and  IsReduced(q)});
\end{verbatim}
The following lines compute the number of points in $P_{0}$ per conics,
the number of conics through points in $P_{0}$, the set $S_{e}$
of number of points in the conics, and then output the results:
\begin{verbatim}
  NumbPoInConics:=[#{p : p in Po| p in q}: q in Conics6]; 
  NumbConicsTroughPo:=[#{q : q in Conics6| p in q}: p in Po]; 
  Se:= Sort(SetToSequence({p : p in NumbPoInConics}));
  return Conics6,NumbConicsTroughPo,Se,NumbPoInConics;
end function;
\end{verbatim}
The second algorithm takes as entry a set $\mathcal{C}$ of irreducible
curves and computes the singularities of the curve $\sum_{q\in\mathcal{C}}q$
with their number and multiplicities. For $r\geq2$, we denote by
$t_{r}$ the number of points of multiplicity $r$ in the curve $\sum_{q\in\mathcal{C}}q$
. 
\begin{verbatim}
function ConfigurationOfCurves(ConfOfCur);
  P2<x,y,z>:=Ambient(ConfOfCur[1]);
  K:=BaseField(P2);
  R1<X>:=PolynomialRing(K);
  AllDeg:=[Degree(q) : q in ConfOfCur];  ma:=Max(AllDeg);
\end{verbatim}
We start by computing the set $SetSchP$ of irreducible components
over the base field $K$ of the intersections of two different curves;
these are $0$-dimensional scheme, maybe not reduced if the intersection
points of the two curves are not transverse:
\begin{verbatim}
  SetSchP:={}; 
  for k1 in [1..#ConfOfCur-1] do
    for k2 in [k1+1..#ConfOfCur] do
      U:=IrreducibleComponents(ConfOfCur[k1] meet ConfOfCur[k2]);
      SetSchP:=SetSchP join {q:q in U};
    end for;
  end for;
\end{verbatim}
We put away the reduced sub-schemes of the $0$-dimensional schemes
in the list $SetSchP$ according to their degree and keep in a list
the non-reduced intersection points:
\begin{verbatim}
  IntPoSch:=[[]: a in [1..ma^2]];
  NotRedIntPo:={};
  for q in SetSchP do
    dq1:=Degree(q);
    if dq1 eq 1 then Append(~IntPoSch[1],q); 
      else
      qr:=ReducedSubscheme(q);
      dq2:=Degree(qr); 
        if dq1 ne dq2 then 
          NotRedIntPo:= {qr} join NotRedIntPo; 
        end if;
      Append(~IntPoSch[dq2],qr);
    end if;
  end for;
  for k in [1..ma^2] do
    W:=IntPoSch[k];
    IntPoSch[k]:=SetToSequence({q : q in W});
  end for;
\end{verbatim}
In the last above lines we removed the eventual repetitions. The $k^{th}$
element in the list $IntPoSch$ is the list of degree $k$ reduced
singular points of the curve $\sum_{q\in\mathcal{C}}q$ over $K$.
For each of these $0$-dimensional schemes $p$, we compute the number
of curves $q$ such that $p$ is a sub-scheme of $q$. In the following
lines, the polynomial $Poly$ is 
\[
Poly=\sum t_{k}X^{k}.
\]
We put a multiplicity $\mu$ equal to the degree of the $0$-dimensional
scheme $p$: if there are $k$ conics containing $p$ and $p_{1},\dots,p_{\mu}$
are the (degree $1$) points in $\bar{K}$ over $p$, each $p_{j}$
has multiplicity $k$.
\begin{verbatim}
  Poly:=R1!0;
  NumCurByPo:=[];
  for mu in [1..ma^2] do
    W:=IntPoSch[mu];
    if #W ne 0 then
      NumCurByPo[mu]:=[#[q : q in ConfOfCur | IsSubscheme(p,q)] : p in W]; 
      Poly:=Poly + &+ [mu*X^u : u in   NumCurByPo[mu]]; 
    end if;
  end for;
\end{verbatim}
Then, we extract the data for the list $T_{r}$ of the integers $t_{2},t_{3},\dots$
\begin{verbatim}
  CoPoly:=Coefficients(Poly);
  Tr:=[[t-1,CoPoly[t]]: t in [1..#CoPoly]| CoPoly[t] ne 0];
\end{verbatim}
Finally, we compute the $H$-constant of the configuration and we
return the datas:
\begin{verbatim}
  Hcst:=((&+AllDeg)^2-&+[q[1]^2*q[2] : q in Tr])/(&+[q[2] : q in Tr]);
  return Tr,Hcst,IntPoSch,NumCurByPo,NotRedIntPo;
end function;
\end{verbatim}

\noindent Xavier Roulleau,
\\Aix-Marseille Universit\'e, CNRS, Centrale Marseille,
\\I2M UMR 7373,  
\\13453 Marseille, France
\\ {\tt Xavier.Roulleau@univ-amu.fr}

http://www.i2m.univ-amu.fr/perso/xavier.roulleau/Site\_Pro/Bienvenue.html
\end{document}